# A Possible Solution for Hilbert's Unsolved 8th Problem: Twin Prime Conjecture

Yuhsin Chen, Yensen Ni[*], and Muyi Chen

## Abstract

We measure whether there are numerous pairs of twin primes (hereafter referred to as twin prime pairs) according to the prime number inferred by sieve of Eratosthenes. In this study, while given a number M= $(6n+5)^2$, we are able to find at least 3 twin prim pairs from the incremental range increased from $(6n+5)^2$ to $[6(n+1)+5]^2$ as n is set from 0 to infinite. Thus, we might prove the Twin Prime Conjecture proposed by de Polignac in 1849. That is, there might have numerous twin prime pairs, indicating that there are numerous prime number p for each natural number k by making p + 2k as the prime number for the case of k = 1.

## Contents

---

[*] Corresponding author: Professor, Department of Management Sciences, Tamkang University, Email: ysni@mail.tku.edu.tw, ysniysni@gmail.com

# 1. Introduction

We measure whether there are numerous pairs of twin primes (hereafter referred to as twin prime pairs) by employing the concept of inferring the prime number proposed by sieve of Eratosthenes. While we set $(6n+5)^2$ as the range for estimating twin prime pairs, we reveal at least three additional twin prime pairs as n (positive integer) is increased by 1. Thus, we confirm the Twin Prime Conjecture made by de Polignac in 1849 by proving that there are numerous prime number p for each natural number k by making p +2k as prime number as k = 1.

In 1990, Hilbert raised 23 mathematical problems at the 1900 International Mathematical Conference held in Paris [2]. The Twin Prime Puzzle (i.e. whether are numerous twin primes existed) is still listed in the $8^{th}$ question, even though Twin Prime Conjecture is proposed by de Polignac in 1849.

To my understanding, these are only three questions unsolved completely among these 23 mathematical problems [2]. As for the recent studies in terms of Hilbert's unsolved $8^{th}$ problems: Twin Prime Conjecture, Y. Zhang proves that bounded gaps between primes are all less than 70 million [7].

After establishing a bound of 70,000,000 on the narrowest gap between primes [7], there are tremendous interests in employing Zhang's argument to lower the bound. T. Tao and dozens of mathematicians worked together to cut the bound down from 70,000,000 to the final value of 246 via an online "Polymath" project as well as more traditional research channels [6]. However, these efforts are based on the conjecture about prime intervals [1]. In other words, the proposition seems changed from whether there are numerous twin prime pairs existed to another proposition: how much the upper limit of the distance between two prime numbers is [3], [5], [6].

Instead of narrowing down the distance between two prime numbers, we focus on the Twin Prime Conjecture: whether there are numerous twin prime pairs. In addition, the question of whether there exist numerous twin prime pairs has been one of the great open questions in number theory for more than a century.

In this study, we endeavor to prove whether there are numerous twin prime pairs. While setting $(6n+5)^2$ as the range[1] for estimating twin prime pairs in accordance

---

[1] According to sieve of Eratosthenes [4], N would be prime number if N would not be divided by the prime number smaller than $\sqrt{N}$. Thus, while setting $\sqrt{N}$ = 6n-1 (prime number) or 6n+1 (prime number), N would be a prime number close to $(6n+5)^2$ (i.e. $(6(n+1)-1)^2$ ), which would not be divided

with Sieve of Eratosthenes [4], we prove that at least three additional twin prime pairs would exist as n (positive integer) is increased by 1. In other words, twin prime pairs would be numerous as n is increased from 1 to $\infty$.

## 2. Associate twin prime (ATP) pairs defined and employed

In this study, we define the pairs (6n-1, 6n+1) as ATP (associate twin prime) pairs whose the first number and the second number in ATP pairs would not be included the multiple of 2 and 3 as shown in Table 1.

For the sake of convenient presentation in this study, we set $ATPn^a$ and $ATPn^b$ as the first number and second number of ATP pairs (6n-1, 6n+1), i.e. 6n-1 and 6n+1 in ATP pairs (6n-1, 6n+1). That is, as for any natural number starting from 4, any multiple of 2 or 3 would be excluded as ATP pairs such as the natural numbers shown in Column 1, Column 3, Column 5, and Column 6 of Table 1. The ATP pairs would be shown in Column 2 and Column 4 of Table 1 since the first number ($ATPn^a$) and the second number ($ATPn^b$) in ATP pairs (ATPn) would be (6n-1, 6n+1).

Table 1.  ATP pairs (6n-1, 6n+1)

| n | Column 1 | Column 2 $ATPn^a$ | Column 3 | Column 4 $ATPn^b$ | Column 5 | Column 6 |
|---|---|---|---|---|---|---|
|  | 6n-2 | 6n-1 | 6n | 6n+1 | 6n+2 | 6n+3 |
| n=1 | 4 | 5 | 6 | 7 | 8 | 9 |
| n=2 | 10 | 11 | 12 | 13 | 14 | 15 |
| n=3 | 16 | 17 | 18 | 19 | 20 | 21 |
| n=4 | 22 | 23 | 24 | 25 | 26 | 27 |
| n=5 | 28 | 29 | 30 | 31 | 32 | 33 |
| n=6 | 34 | 35 | 36 | 37 | 38 | 39 |
| n=7 | 40 | 41 | 42 | 43 | 44 | 45 |
| n=8 | 46 | 47 | 48 | 49 | 50 | 51 |
| n=9 | 52 | 53 | 54 | 55 | 56 | 57 |
| n=10 | 58 | 59 | 60 | 61 | 62 | 63 |
| n=11 | 64 | 65 | 66 | 67 | 68 | 69 |
| n=12 | 70 | 71 | 72 | 73 | 74 | 75 |
| n=13 | 76 | 77 | 78 | 79 | 80 | 81 |
| n=14 | 82 | 83 | 84 | 85 | 86 | 87 |
| n=15 | 88 | 89 | 90 | 91 | 92 | 93 |
| n=16 | 94 | 95 | 96 | 97 | 98 | 99 |
| n=17 | 100 | 101 | 102 | 103 | 104 | 105 |
| ⋮ | ⋮ | ⋮ | ⋮ | ⋮ | ⋮ | ⋮ |

by the prime number small than $\sqrt{N}$ such as 6n+1 due to that 6n+1 is smaller than 6n+5. Thus, we set $(6n+5)^2$ as the range for estimating Twin Prime pairs (pairs of twin primes) in accordance with Sieve of Eratosthenes

Thus, if a natural number would not be divided by 2 or 3, then the natural number would be presented as 6n±1. As a result, prime numbers except 2 and 3 would be presented as 6n±1, where n is a positive integer. As a result, the difference between 6n-1 and 6n+1 is exactly 2, i.e., the pairs (6n+1 and 6n-1) for n=1 to ∞ would be regarded as ATP pairs.

In addition, the ATP pairs (6n-1, 6n+1) would contain all twin prime pairs except the pair (3, 5) such as twin prime pair (5, 7) as n =1, twin prime pair (11, 13) as n=2, twin prime (17, 19) as n =3, twin prime pair (29, 31) as n=5, twin prime pair (41, 43) as n=7, twin prime pair (59, 61) as n=10. However, the ATP pairs (6n-1, 6n+1) are not all twin prime pairs for any n (positive integer). For example, the ATP pairs (6n-1, 6n+1) are not twin-prime pairs for the pair (23, 25) as n=4, the pair (35, 37) as n =6, the pair (47, 49) as n=8, and the pair (53, 55) as n=9. That is, as for the ATP pairs (6n-1, 6n+1) for n is from 1 to 10, we show that twin prime pairs belong to ATP pairs except for the twin prime pair (3, 5); however, ATP pairs are not all twin prime pairs. In other words, if and only if $ATPn^a$ and $ATPn^b$ are all prime numbers, then the ATP pair (ATPn) would be regarded as a twin prime pair.

The above concern is to explore whether we are able to find the evidence of numerous twin prime pairs due to that the number of twin prime pair would be increased continuously as the selected incremental range expanded continuously. Thus, we would explain the selection for the incremental range in the following section.

### 3. Selection for Incremental Range

By calculating the number of ATP pairs within a certain range, we are able to explain why the number of ATP pairs would increase as the range is increased incrementally (i.e. the range is increased incrementally from $(6n+5)^2$ to $(6(n+1)+5)^2$ for n=1 to ∞). That is, we argue that the number of twin prime pairs would be increased infinitely when selecting the range is expanded infinitely. Then, the conjecture of the twin prime would be confirmed in this study.

According to sieve of Eratosthenes, N would be confirmed as prime number if N would not be divided by the prime number not bigger than $\sqrt{N}$. Thus, while setting $\sqrt{N}$ = 6n-1 (prime number) or 6n+1 (prime number), N would be a prime number close to $(6n+5)^2$ (i.e. $(6(n+1)-1)^2$ = N), which would not be divided by any prime number small than $\sqrt{N}$ such as 6n+1 due to that 6n+1 is smaller than 6n+5 (i.e. 6n+1 < 6n+5 = $\sqrt{N}$ ).

As mentioned above, if and only if ATPn[a] and ATPn[b] are all prime numbers, then the ATP pair (ATPn) would be regarded as a twin prime pair. Thus, the Twin Prime pairs might be estimated as long as we could examine whether all these ATPn[a] and ATPn[b] less than (6n+5)$^2$ (i.e. (ATP(n+1)[a])$^2$ = (6(n+1)-1)$^2$) would not be divided by any prime number less than ATP(n+1)[a] such as 6n+1 less than 6n+5 (i.e. 6(n+1)-1). Thus, the following context is set (6n+5)$^2$ as the selection for the incremental range in this study, which would be convenient for us to explore whether there are unlimited twin prime pairs in this study.

While estimating the twin prime pairs for the range smaller than (6n+5)$^2$ in accordance with sieve of Eratosthenes, we calculate the number of ATP pairs for the range of (6n+5)$^2$ for n as n =1 to unlimited (hereafter referred to as ATPGn). Due to ATP pairs would be calculated as a cycle of every six integers, the ATPGn would include $\left\{\frac{[6n+5]^2-1}{6} - 1\right\}$ ATP pairs. That is, ATPG1 (i.e. ATPGn as n=1) would have 19 ATP pairs as the range is less than (6n+5)$^2$= 121 (i.e. n=1). Similarly, ATPG2=47 as the range is less than (6n+5)$^2$= 289 (i.e. n=2), ATPG3=87 as the range is less than (6n+5)$^2$= 529 (i.e. n=3), ATPG4=139 as the range is less than (6n+5)$^2$= 841 (i.e. n=4), and so on.

Thus, if the twin prime numbers for ATPGn would be increased more than 1 as long as n is incremented by 1 each time, then we are able to prove Twin Prime Conjecture, i.e. (i.e. there are numerous twin primes pairs existed).

We then calculate the increasing rate of ATPGn (i.e. ATPGn%) as n is incremented by 1, and the formula is shown as below.

ATPGn%=$\left\{\frac{[ATP(n+1)^a]^2-1}{6} - 1\right\} / \left\{\frac{[ATP(n)^a]^2-1}{6} - 1\right\}$

=$\left\{\frac{[6(n+1)-1]^2-1}{6} - 1\right\} / \left\{\frac{[6(n)-1]^2-1}{6} - 1\right\}$

=$\left\{\frac{[6n+5]^2-1}{6} - 1\right\} / \left\{\frac{[6n-1]^2-1}{6} - 1\right\}$

=$\left\{\frac{36n^2+60n+25-1}{6} - 1\right\} / \left\{\frac{36n^2-12n+1-1}{6} - 1\right\}$

$$=\frac{6n^2+10n+3}{6n^2-2n-1}$$

As mentioned above in section 2, the first number and the second number of an ATP pair defined in this study would be started from the integer 5 (i.e. the smallest first number of ATP pair (6n-1, 6n+1) would be 5). As a result, the Associate Twin Prime Group (ATPG0) would be equal to 3 for n=0 (i.e. the ATP pairs (5, 7), (11, 13), and (17, 19)) as the range smaller than $(6n+5)^2 = 25$ as n=0.

As a result, the ATPGn would be presented as ATPGn-1· ATPGn%, and expand to ATPG0 · ATPG1% · ATPG2% · ATPG3% · .... ATPGn%, which could be calculated as shown below.

$$ATPGn = 3 \cdot \prod_{i=1}^{n} \frac{6i^2 + 10i + 3}{6i^2 - 2i - 1}$$

where the initial value = 3 as n=0 (i.e $(6n+5)^2 = 25$, if n=0).

By employing the above equation, we would get the ATPG0=3 (i.e. the initial value =3), ATPG1= ATPG0 · ATPG1% = $3 \cdot \frac{19}{3}$ = 19, ATPG2 = ATPG1 · ATPG2% = $3 \cdot \frac{19}{3} \cdot \frac{47}{19}$ = 47, and so on.

In addition, we would further measure the equation: $TPEn = ATPGn \cdot TPRn$, where TPEn is the estimated number of twin prime pairs accumulated to the range less than $(6n+5)^2$ for n = 1 to unlimited, ATPGn is the number of ATP pairs as integer number less than $(6n+5)^2$ for n=1 to unlimited, and TPRn is the probability of ATP pairs being twin prime pairs within the range $(6n+5)^2$ for n=1 to unlimited. As a result, we are able to measure the TPEn (The estimated number of Twin Prime pairs accumulated from the range $(6n+5)^2$ as n=0 (initial value =3) to the range $(6n+5)^2$ as n= unlimited.

In this study, we define TPEn as the estimated number of twin prime pairs accumulated to the range less than $(6n+5)^2$ for n = 1 to unlimited, ATPGn as the number of ATP pairs as integer number less than $(6n+5)^2$ for n=1 to unlimited, and TPRn as the probability of ATP pairs being twin prime pairs within the range $(6n+5)^2$ for n=1 to unlimited. Then, by measuring the equation: $TPEn = ATPGn \cdot TPRn$, we are able to measure the TPEn from the range $(6n+5)^2$ as n=1 [2] to the range $(6n+5)^2$ as n=

---

[2] In fact, the initial value =3 for TPEn as n=0, since there are three twin prime pairs, (3, 5), (11, 13), and (17, 19) less than 25, (i.e. $(6n+5)^2$ as n=0).

unlimited. In other words, we would measure whether the twin prime pairs would be increased unlimitedly as the incremental range, $(6n+5)^2$ increased from 1 to limited in this study. Thus, we would measure the probability of ATP pairs being twin prime pairs within the range $(6n+5)^2$ for n=1 to unlimited in the following section.

## 4. The probability of ATP pairs being Twin prime pairs

In addition, ATPn would include a pair of natural numbers (ATPn$^a$, ATPn$^b$). While selecting any natural number x in ATPn (such as ATPn$^a$ = 5 as n=1), two ATP pairs out of five ATP pairs would not be twin prime pairs as shown in Figure 1; similarly, while selecting any natural number x in ATPn (such as ATPn$^b$ = 7 as n=1), two ATP pairs out of seven ATP pairs would not be twin prime pairs as shown in Figure 1. In other words, the probability of excluding twin prime pairs from ATP pairs would be 2/x. For example, the probability of excluding twin prime pairs from ATP pairs would be 2/5 as x is 5 (i.e. ATP1$^a$=5) as well as 2/7 as x is 7 (i.e. ATP1$^b$=7) as shown in Figure 1.

| ATP Pairs | A | b | | ATP Pairs | a | b | |
|---|---|---|---|---|---|---|---|
| 1 | 5 | 7 | | 1 | 5 | 7 | |
| 2 | 11 | 13 | | 2 | 11 | 13 | |
| 3 | 17 | 19 | | 3 | 17 | 19 | |
| 4 | 23 | 25 | 5 pairs | 4 | 23 | 25 | |
| 5 | 29 | 31 | | 5 | 29 | 31 | 7 pairs |
| 6 | 35 | 37 | | 6 | 35 | 37 | |
| 7 | 41 | 43 | | 7 | 41 | 43 | |
| 8 | 47 | 49 | | 8 | 47 | 49 | |
| 9 | 53 | 55 | 5 pairs | 9 | 53 | 55 | |
| 10 | 59 | 61 | | 10 | 59 | 61 | |
| 11 | 65 | 67 | | 11 | 65 | 67 | |
| 12 | 71 | 73 | | 12 | 71 | 73 | 7 pairs |
| 13 | 77 | 79 | | 13 | 77 | 79 | |
| 14 | 83 | 85 | 5 pairs | 14 | 83 | 85 | |
| 15 | 89 | 91 | | 15 | 89 | 91 | |
| 16 | 95 | 97 | | 16 | 95 | 97 | |
| 17 | 101 | 103 | | 17 | 101 | 103 | |
| 18 | 107 | 109 | | 18 | 107 | 109 | |
| 19 | 113 | 115 | 5 pairs | 19 | 113 | 115 | 7 pairs |
| 20 | 119 | 121 | | 20 | 119 | 121 | |
| 21 | 125 | 127 | | 21 | 125 | 127 | |
| 22 | 131 | 133 | | 22 | 131 | 133 | |
| 23 | 137 | 139 | | 23 | 137 | 139 | |
| 24 | 143 | 145 | 5 pairs | 24 | 143 | 145 | |
| 25 | 149 | 151 | | 25 | 149 | 151 | |
| 26 | 155 | 157 | | 26 | 155 | 157 | 7 pairs |
| 27 | 161 | 163 | | 27 | 161 | 163 | |
| 28 | 167 | 169 | | 28 | 167 | 169 | |
| 29 | 173 | 175 | 5 pairs | 29 | 173 | 175 | |
| 30 | 179 | 181 | | 30 | 179 | 181 | |
| 31 | 185 | 187 | | 31 | 185 | 187 | |
| 32 | 191 | 193 | | 32 | 191 | 193 | |
| 33 | 197 | 199 | | 33 | 197 | 199 | 7 pairs |
| 34 | 203 | 205 | 5 pairs | 34 | 203 | 205 | |
| 35 | 209 | 211 | | 35 | 209 | 211 | |

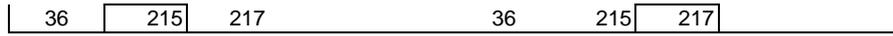

Figure 1. Excluding Twin Prim pairs from ATP pairs for Prime numbers 5 and 7

In Figure 1, the numbers in black boxes are excluded as prime numbers due to these numbers would be divided by 5 or 7. As either the first number or second number in an ATP pair (6n-1, 6n+1) is not prime numbers, this ATP pair would not be a Twin Prime pair. We also find that excluding twin prims pairs from ATP pairs would be rather regular, i.e. the probability of excluding twin prime pair from ATP pairs would be 2/x.

According to the screening approach proposed by sieve of Eratosthenes, some ATP pairs would be excluded from twin prime pairs, since the multiples of prime numbers smaller than $\sqrt{N}$ are shown in the first numbers or second numbers in ATP pairs.

In fact, some of the first number and second number in ATP are not prime numbers instead of composite numbers, which would not be included as twin prime pair. For example, the number 35 in ATP pair would be excluded as well as the numbers 175 and 245 in ATP pairs would be excluded as well due to that these numbers are composite numbers.

Thus, we could prove "at least" a twin prime pair existed for the incremental range from $(6n+5)^2$ to $(6(n+1)+5)^2$ for n =1, 2, 3,….infinite, since initial value =3 as n=0 [3]. In addition, we argue that our proof would be more persuasive if we are able to prove "at least" one twin prime pair existed within the incremented range. In addition, we document that our estimation would be underestimated, which would be explained as below:

As mentioned above, while setting $\sqrt{N}$ = 6n-1 or 6n+1, N would be a prime number close to $(6n+5)^2$ (i.e. $(6(n+1)-1)^2$ = N), which would not be divided by any prime number small than $\sqrt{N}$ such as 6n+1 due to that 6n+1 is smaller than 6n+5 (i.e. 6n+1 < 6n+5 = $\sqrt{N}$). As a result, ATPG1 (i.e. ATPGn as n=1) would not be divided by the prime number small than $\sqrt{121}$ (i.e. 5 & 7), ATPG2 (i.e. ATPGn as n=2) would not be divided by the prime number small than $\sqrt{289}$ (i.e. 5, 7, 11, and 13),

---

[3] The Associate Twin Prime Group (ATPG) would be equal to 3 including (5, 7), (11, 13), and (17, 19), for the integer smaller than $5^2$ (i.e. $(6n+5)^2$ = 25 as n=0). In addition, Twin Prime pairs would be equal to 3 as well, since (5, 7), (11, 13), and (17, 19) are all twin prime pairs.

ATPG3 (i.e. ATPGn as n=3) would not be divided by the prime number small than $\sqrt{529}$ (i.e. 5, 7, 11, 13, 17, and 19), ATPG4 (i.e. ATPGn as n=4) would not be divided by the prime number small than $\sqrt{841}$ (i.e. 5, 7, 11, 13, 17, 19, 23), ATPG5 (i.e. ATPGn as n=5) would not be divided by the prime number small than $\sqrt{1225}$ (i.e. 5, 7, 11, 13, 17, 19, 23, 29, 31), and so on. The above explanation would be beneficial for us to measure the probability of Cases I – IV, as shown below. In addition, the probability of Case I would be regarded as a conservative probability as illustrated below as well.

In this study, the probability of Case I: ATPn$^a$ and ATPn$^b$ are all prime numbers would be $TPRn^{ab}$ (The ratio of twin prime pairs over ATP pairs; Twin Prime Ratio) could be measured as shown below.

$$TPRn^{ab} = (1 - \frac{2}{\text{ATPn}^a})(1 - \frac{2}{\text{ATPn}^b})$$ for Case I: ATPn$^a$ and ATPn$^b$ are all prime numbers.

We also find that excluding twin prims pairs from ATP pairs would be rather regular, i.e. the probability of excluding twin prime pair from ATP pairs would be 2/x. In addition, x would be the first number or second number in ATP pair (ATPn$^a$, ATPn$^b$).

However, we argue that there are four cases for ATP pair (ATPn$^a$, ATPn$^b$) including Case I, (prime number, prime number), Case II, (prime number, composite number), Case III, (composite number, prime number), and Case IV (composite number, composite number). Thus, we are able to measure various probabilities for Cases II – IV as shown below.

$$TPRn^a = (1 - \frac{2}{\text{ATPn}^a})$$ for Case II: ATPn$^a$ is prime number but ATPn$^b$ is not prime numbers.

$$TPRn^b = (1 - \frac{2}{\text{ATPn}^b})$$ for Case III: ATPn$^a$ is not prime number but ATPn$^b$ is prime numbers.

$TPRn^0 = 1$ for Case IV: ATPn$^a$ and ATPn$^b$ are all not prime numbers.

In addition, we are able to compare the probability these cases as shown below.

$$\left(1 - \frac{2}{\text{ATPn}^a}\right)\left(1 - \frac{2}{\text{ATPn}^b}\right) < 1 - \frac{2}{\text{ATPn}^a} < 1 - \frac{2}{\text{ATPn}^b} < 1$$

That is, $TPRn^{ab} < TPRn^a < TPRn^b < TPRn^0$

Moreover, we employ the probability, $TPRn^{ab}$, rather than $TPRn^a$, $TPRn^b$, and $TPRn^0$ to measure whether there are numerous twin prime pairs existed, so our estimator would be regarded as conservative estimation.

$$TPRn = \prod_{i=1}^{n}(1 - \frac{2}{ATPi^a})(1 - \frac{2}{ATPi^b})$$

$$= \prod_{i=1}^{n}(1 - \frac{2}{6i-1})(1 - \frac{2}{6i+1})$$

$$= \prod_{i=1}^{n}(\frac{6i-1-2}{6i-1})(\frac{6i+1-2}{6i+1})$$

$$= \prod_{i=1}^{n}\frac{(6i-3)(6i-1)}{(6i-1)(6i+1)}$$

$$= \prod_{i=1}^{n}\frac{(6i-3)}{(6i+1)}$$

### 5. Evidence for the existence of infinite Twin Prime Pairs

In summary, the conservative estimation of Twin Prime pairs could be measured as that the number of ATP Group (ATPGn) multiplied by the probability of Twin Prime pairs over ATP pairs (TRPn).

Thus, TPEn = ATPGn・TPRn , where TPEn is the estimated number of twin prime pairs accumulated to the range less than $(6n+5)^2$ for n = 1 to unlimited, ATPGn is the number of ATP pairs as integer number less than $(6n+5)^2$ for n=1 to unlimited, and TPRn is the probability of ATP pairs being twin prime pairs within the range $(6n+5)^2$ for n=1 to unlimited. We then employ $TPRn^{ab}$ instead of $TPRn^a$, $TPRn^b$, and $TPRn^0$ to measure the above equation, and the numbers of Twin Prime pairs estimated would be rather conservative as well.

Due to $ATPGn = 3 \cdot \prod_{i=1}^{n} \frac{6i^2+10i+3}{6i^2-2i-1}$ and $TPRn^{ab} = \prod_{i=1}^{n} \frac{(6i-3)}{(6i+1)}$,

TPEn

$$= 3 \cdot \prod_{i=1}^{n} (\frac{6i^2+10i+3}{6i^2-2i-1})(\frac{6i-3}{6i+1})$$

$$= 3 \cdot \prod_{i=1}^{n} (\frac{36i^3+60i^2+18i-18i^2-30i-9}{36i^3-12i^2-6i+6i^2-2i-1})$$

$$= 3 \cdot \prod_{i=1}^{n} (\frac{36i^3+42i^2-12i-9}{36i^3-6i^2-8i-1})$$

In order to prove Twin Prime Conjecture, the above math equation would be simplified as below.

$$= 3 \times \prod_{i=1}^{n} (\frac{36i^3+42i^2-12i-9}{36i^3-6i^2-8i-1})$$

$$= 3 \times \prod_{i=1}^{n} [\frac{(36i^3-6i^2-8i-1)+(48i^2-4i-8)}{36i^3-6i^2-8i-1}]$$

$$= 3 \times \prod_{i=1}^{n} [1 + \frac{((48i^2-4i-8)}{36i^3-6i^2-8i-1}]$$

$$= 3 \times \prod_{i=1}^{n} [1 + \frac{(36i^2-6i-8-\frac{1}{i})}{36i^3-6i^2-8i-1} + \frac{(12i^2+2i+\frac{1}{i})}{36i^3-6i^2-8i-1}]$$

$$= 3 \times \prod_{i=1}^{n} [1 + \frac{(36i^2-6i-8-\frac{1}{i})}{i(36i^2-6i-8-\frac{1}{i})} + \frac{(12i^2+2i+\frac{1}{i})}{36i^3-6i^2-8i-1}]$$

$$= 3 \times \prod_{i=1}^{n} [1 + \frac{1}{i} + \frac{(12i^2+2i+\frac{1}{i})}{36i^3-6i^2-8i-1}]$$

$$> 3 \times \prod_{i=1}^{n} [1 + \frac{1}{i}] \qquad \because i \geq 1, \quad \frac{(12i^2+2i+\frac{1}{i})}{36i^3-6i^2-8i-1} > 0$$

Due to $\frac{(12i^2+2i+\frac{1}{i})}{36i^3-6i^2-8i-1}$ would not be negative, the Twin Prime pairs would be greater than $3 \times \prod_{i=1}^{n}[1+\frac{1}{i}]$. After simplifying $3 \times \prod_{i=1}^{n}[1+\frac{1}{i}]$, we would get $3 \times \prod_{i=1}^{n}[1+\frac{1}{i}]$ would be equal to 3n +3 as shown below.

$$3 \times \prod_{i=1}^{n}[1+\frac{1}{i}] = 3 \times \prod_{i=1}^{n}[\frac{i+1}{i}] = 3 \times \frac{2}{1} \times \frac{3}{2} \times \ldots \times \frac{n}{n-1} \times \frac{n+1}{n}$$

$= 3 (n+1) = 3n+3$

After we differentiate the above equation, we would determine that this conservative estimator is an incremental function

$\frac{d}{dn}(3n+3)=3$

The above inference indicates that three Twin Prime pairs would be increased as n is increased by 1.

Thus, we then prove the Twin Prime Conjecture proposed by de Polignac in 1849. That is, there are numerous Twin Prime pairs, indicating that there are numerous prime number p for each natural number k by making p +2k as prime number as well for the case of k = 1.

In addition, we measure TPEn (i.e. TPEn = ATPGn・TPRn, while setting TPRn = $TPRn^{ab}$), TPAn (i.e. the actual twin prime pairs), and TPEn-S (i.e. the more conservative estimation for measuring TPE due to excluding non-negative term as well as setting TRRn = $TPRn^{ab}$) for n=1 to 50 as shown in Table 2. In addition, we observe that TPEn is less than TPAn as well as TPEn-S is far less than TPAn. Even so, we confirm that there are numerous twin prime pairs, as the results shown for TPEn-S and TPEn.

Table 2. The measurement for TPEn, TPAn, and TPE-S

|  |  | ATPGn | TPRn | TPEn | TPAn | TPEn-S |
| --- | --- | --- | --- | --- | --- | --- |

| n | Incremental range $(6n+5)^2$ | The number of ATP pairs $\left\{\dfrac{[6n+5]^2-1}{6}-1\right\}$ | Conservation (1) The conservative probabilities of $TPRn^{ab}$: $\prod_{i=1}^{n}(1-\dfrac{2}{ATPi^a})(1-\dfrac{2}{ATPi^b})$ | TPEn | Actual TP numbers | Conservative & Simplified (1) The conservative probabilities of $TPRn^{ab}$: $\prod_{i=1}^{n}(1-\dfrac{2}{ATPi^a})(1-\dfrac{2}{ATPi^b})$ (2) Excluding the non-negative term |
|---|---|---|---|---|---|---|
| 1 | 121 | 19 | 42.86% | 8 | 9 | 6 |
| 2 | 289 | 47 | 29.67% | 14 | 18 | 9 |
| 3 | 529 | 87 | 23.42% | 20 | 24 | 12 |
| 4 | 841 | 139 | 19.68% | 27 | 32 | 15 |
| 5 | 1,225 | 203 | 17.14% | 35 | 40 | 18 |
| 6 | 1,681 | 279 | 15.28% | 43 | 52 | 21 |
| 7 | 2,209 | 367 | 13.86% | 51 | 66 | 24 |
| 8 | 2,809 | 467 | 12.73% | 59 | 79 | 27 |
| 9 | 3,481 | 579 | 11.81% | 68 | 92 | 30 |
| 10 | 4,225 | 703 | 11.03% | 78 | 109 | 33 |
| 11 | 5,041 | 839 | 10.37% | 87 | 127 | 36 |
| 12 | 5,929 | 987 | 9.80% | 97 | 142 | 39 |
| 13 | 6,889 | 1,147 | 9.31% | 107 | 159 | 42 |
| 14 | 7,921 | 1,319 | 8.87% | 117 | 173 | 45 |
| 15 | 9,025 | 1,503 | 8.48% | 127 | 190 | 48 |
| 16 | 10,201 | 1,699 | 8.13% | 138 | 209 | 51 |
| 17 | 11,449 | 1,907 | 7.81% | 149 | 226 | 54 |
| 18 | 12,769 | 2,127 | 7.53% | 160 | 243 | 57 |
| 19 | 14,161 | 2,359 | 7.27% | 171 | 262 | 60 |
| 20 | 15,625 | 2,603 | 7.03% | 183 | 277 | 63 |
| 21 | 17,161 | 2,859 | 6.80% | 195 | 297 | 66 |
| 22 | 18,769 | 3,127 | 6.60% | 206 | 324 | 69 |
| 23 | 20,449 | 3,407 | 6.41% | 218 | 346 | 72 |
| 24 | 22,201 | 3,699 | 6.23% | 231 | 375 | 75 |
| 25 | 24,025 | 4,003 | 6.07% | 243 | 401 | 78 |
| 26 | 25,921 | 4,319 | 5.91% | 255 | 417 | 81 |
| 27 | 27,889 | 4,647 | 5.77% | 268 | 443 | 84 |
| 28 | 29,929 | 4,987 | 5.63% | 281 | 466 | 87 |
| 29 | 32,041 | 5,339 | 5.50% | 294 | 491 | 90 |
| 30 | 34,225 | 5,703 | 5.38% | 307 | 526 | 93 |
| 31 | 36,481 | 6,079 | 5.27% | 320 | 553 | 96 |
| 32 | 38,809 | 6,467 | 5.16% | 334 | 580 | 99 |
| 33 | 41,209 | 6,867 | 5.05% | 347 | 601 | 102 |
| 34 | 43,681 | 7,279 | 4.95% | 361 | 629 | 105 |
| 35 | 46,225 | 7,703 | 4.86% | 374 | 654 | 108 |
| 36 | 48,841 | 8,139 | 4.77% | 388 | 684 | 111 |
| 37 | 51,529 | 8,587 | 4.69% | 402 | 722 | 114 |
| 38 | 54,289 | 9,047 | 4.60% | 417 | 747 | 117 |
| 39 | 57,121 | 9,519 | 4.53% | 431 | 779 | 120 |
| 40 | 60,025 | 10,003 | 4.45% | 445 | 810 | 123 |
| 41 | 63,001 | 10,499 | 4.38% | 460 | 834 | 126 |
| 42 | 66,049 | 11,007 | 4.31% | 474 | 867 | 129 |
| 43 | 69,169 | 11,527 | 4.24% | 489 | 894 | 132 |
| 44 | 72,361 | 12,059 | 4.18% | 504 | 937 | 135 |
| 45 | 75,625 | 12,603 | 4.12% | 519 | 964 | 138 |
| 46 | 78,961 | 13,159 | 4.06% | 534 | 994 | 141 |
| 47 | 82,369 | 13,727 | 4.00% | 549 | 1,036 | 144 |
| 48 | 85,849 | 14,307 | 3.94% | 564 | 1,074 | 147 |
| 49 | 89,401 | 14,899 | 3.89% | 580 | 1,108 | 150 |
| 50 | 93,025 | 15,503 | 3.84% | 595 | 1,147 | 153 |

Note: TPEn = ATPGn・TPRn, while setting TPRn = $TPRn^{ab}$ for n=1 to 50

TPAn is the actual twin prime pair for n=1 to 50

TPEn-S is the more conservative estimation for measuring TPE due to excluding non-negative term as well as setting TRRn = $TPRn^{ab}$ for n=1 to 50.

In addition, TPAn (Actual TP numbers) might have the following two concerns. One is the probabilities measured by $TPRn^a, TPRn^b, TPRn^0, or TPRn^{ab}$ might be depended on various cases shown in (ATPn$^a$, ATPn$^b$), and the other is to include the non-negative term as mentioned above.

In addition, we plot the trend for TPEn, TPAn, and TPEn-S. In fact, while setting $(6n+5)^2$ as the incremental range for measuring the numbers of twin prime pairs, we observe that increasing convex curves are shown for TPEn and TPAn, and increasing linear curve is shown for TPn-S due to that TPE-S=3n+3. Moreover, TPAn is much higher than TPEn-S and TPEn. While n is approaching unlimited, we are able to confirm numerous twin prim pars for TPEn-S. However, the numbers of twin prime pairs measured by TPEn-S are far less than those measured by TPEn. Likewise, the numbers of twin prime pairs measured by TPEn are far less than those measured by TPAn, indicating that confirming numerous twin prime pairs existed is rather persuasive in this study.

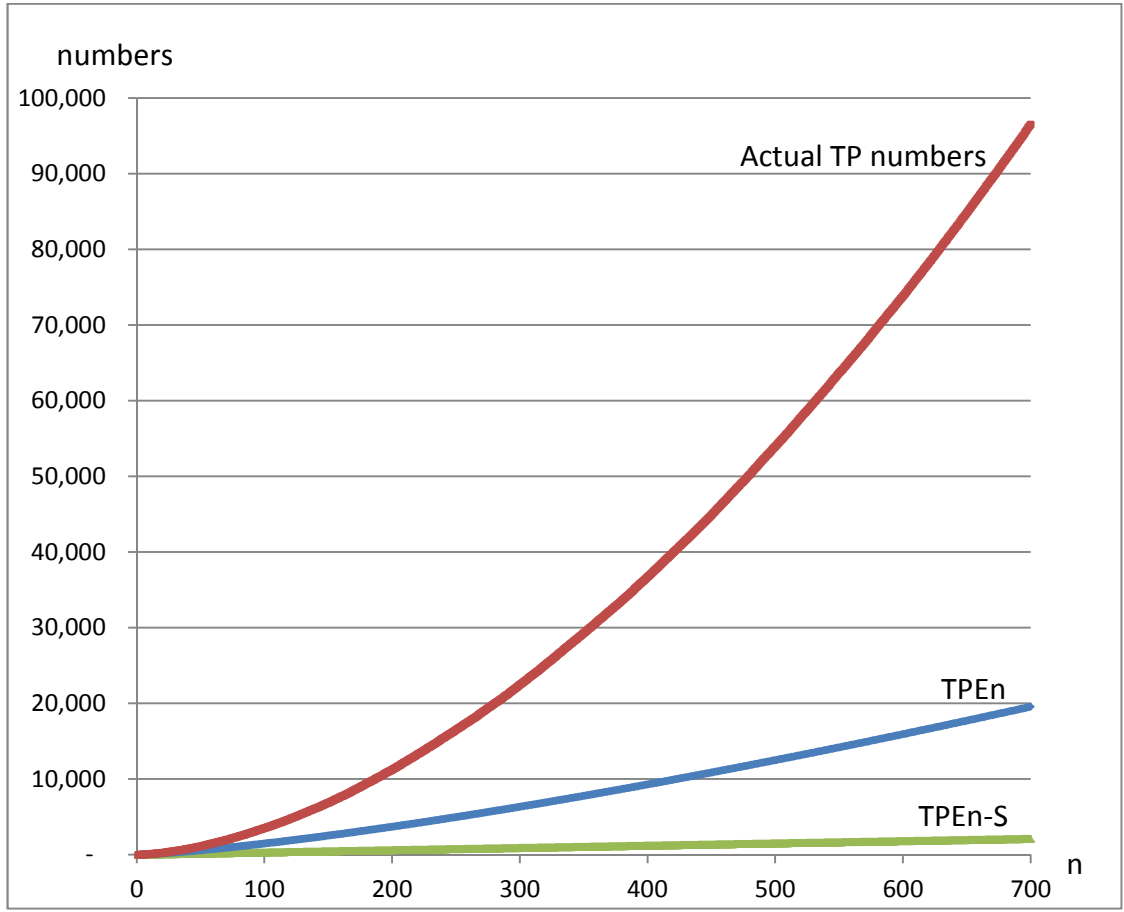

Figure 2. The trend for TPEn, TPAn, and TPEn-S

Note: TPAn: actual twin prime numbers by taking the following two conditions into account
TPEn: conservative twin numbers measuring by the equation ATPGn・TPRn = TPEn
Where TPRn is set as $\prod_{i=1}^{n}(1 - \frac{2}{ATPi^a})(1 - \frac{2}{ATPi^b})$
TPEn-S: very conservative and simplified measurement for twin primes as shown above. (1)
TPRn is set as $\prod_{i=1}^{n}(1 - \frac{2}{ATPi^a})(1 - \frac{2}{ATPi^b})$  (2) Excluding the non-negative term

In fact, the above estimation would be regarded as a rather conservative estimation, and this conservative estimation mainly comes from two parts. One is that we employ $TPRn^{ab}$ instead of TPRn$^a$, TPRn$^b$, and TPRn$^0$ do to measure TPEn-S. The other is that we simplify mathematical presentation by excluding the non-negative term $\frac{(12i^2+2i+\frac{1}{i})}{36i^3-6i^2-8i-1}$ in the mathematical equation. Even so, we still prove that there are numerous prime number p for each natural number k by making p +2k as prime number as well for the case of k = 1.

In addition, if we employ $TPRn^{ab}$, TPRn$^a$, TPRn$^b$, or TPRn$^0$ depended on the actual condition occurred as well as include a non-negative term in the mathematical equation, the numbers of Twin Prime pairs would be increased than the conservative estimation measured above as well as the result would be closer to the actual result. As a result, how to get the estimated result close to the actual result would be worthwhile for further studies.

## 6. Conclusion

We extend the sieve of Eratosthenes to calculate whether there are numerous numbers of Twin Prime pairs. In order to make our revealed result reliable even persuaded, we use very conservative estimation.

In this study, we prove that there are at least 3 additional twins prime pairs increased while setting $(6n+5)^2$ as the range of estimating the additional increase of twins prime pairs as n is increased by 1.

Furthermore, although estimating the twin prime pairs conservatively, we still prove that there are numerous twin prime pairs, since there are at least 3 additional twin prime pairs increased while setting $(6n+5)^2$ as the range of estimating the additional increase of twins prime pairs as n is increased by 1 for n =1, 2, ………..∞.

We then confirm that additional twin prime pairs would be increased as n is increased as well. Due to that n would be increased from 1 to ∞, twin prime pairs could be increased to infinite as well. That it, Hilbert's unsolved 8$^{th}$ problems: Twin Prime Conjecture has been resolved in this study.

In fact, Twin Prime Conjecture solved might be based on two main concerns. One is to delete the non-twin prime pairs over ATP pairs would be regular (i.e. 2/x shown in context). The other is the incremental range from $(6n+5)^2$ to $(6(n+1)+5)^2$ for n =1, 2, 3,….infinite (i.e. the incremental range increased could be regarded as regular as well). We argue that Twin Prime Conjecture solved might be due to these two concerns are taken into account deliberately in the study.

Due to that, the estimation of whether there are infinite numbers of Twin Prim pairs is rather conservative and even simplified, the more accurate estimation would be taken into account for further studies.